# Not served on a silver platter!
## Access to online mathematics information in Africa

Anders Wändahl
"e-Math for Africa"
anders@golonka.se

**Abstract**

*This paper argues that, contrary to the beliefs of many, the amount of mathematics information available for African researchers, including electronic scientific journals and databases, is indeed substantial. However, whereas information resources are served on a silver platter at the "northern" universities, researchers at universities in developing countries have to work hard for their treats.*

*The scientific information available for low-income countries is scattered among a large number of providers, websites, access methods, price models, and country- or institution-specific programmes. It's therefore quite hard for individual researchers to see the whole picture and establish what actually is available and what is not.*

*This paper presents a number of key information sources for mathematics. Some are aimed primarily at disciplines other than mathematics but nevertheless contain extremely important and high-ranking mathematics journals. Brief instructions are also given on how to register for and maintain access to various relevant and useful resources.*

## Why is scientific information important?

Access to scientific information is not only an undisputable condition for being able to perform good research, but also a means to participate generally in the global scientific community. One of the contributing factors for brain-drain from African institutions is the scarcity of library resources such as books, journals and databases, as was pointed out in a report from the African Network of Scientific and Technological Institutions (ANSTI) 2005:

"The availability of journals also has profound influence on the motivation of staff and their career progression. The absence of journals and good library facilities means inability to do research, publish and contribute to knowledge. The staff will not be aware of current trends in his/her field and over some time his/ her knowledge becomes obsolete and career prospects dim. Since most scientists will not be able to fund their own research or journals, there is the tendency to leave the institution once they find that the library facilities are not very good. ***Thus, the observation made in this study on the state of libraries, may also be one of the contributing factors to the brain-drain from African institutions***."[1]

---

[1] ANSTI (2005). *State of Science and Technology Training Institutions in Africa* [Online], p. 39
Available at: http://www.ansti.org/reports/state%20of%20science%20training%20in%20africa.pdf
[Accessed Mar 28th, 2009].

Although the inadequacy of library facilities indeed is a problem, the situation is not as bad as it seems. One of the recent dominant movements on the publishing scene is Open Access, provision of free access for everyone. Additionally, many different programmes have been set up recently to improve access to scientific information for developing countries. These programmes focus on subsidized or free access to journals that readers normally must subscribe to. Only a few are aimed at mathematics information specifically, but most include some mathematics resources.

Generally speaking, the scientific information available for developing countries is scattered among a large number of web sites, access methods, price models, country or institutional specific deals, and it's quite hard for the single developing-country researcher to see the whole picture and establish what is actually available and what is not.

In the "north", university libraries can allocate two or three fulltime librarians to the task of cataloguing, classifying, organizing, and in other ways packaging all the available information resources in order to make the picture complete and unambiguous. One part of the package is usual a special database that lists and gives access to the full-text of available journals. This is what I mean with "served on a silver platter" in the title of this paper. Very few developing country universities have enough money and/or human resources to maintain such a journal database.

What makes mathematics a little different from other disciplines is that the literature tends to last longer. A journal paper in cutting-edge biosciences may be 'old' already at the time of publication. Papers published in the 1970s and 1980s are in many cases considered outdated and may even be incorrect in the light of recent research. In mathematics the situation is different. Classic mathematics books remain classics for decades, and journal articles remain valid as long as the mathematical truths were not false in the first place. Mathematicians rely very much on old material and can now and then find "hidden pearls", that is, old methods that are applicable to new problems. [2]

## The role of reference databases in finding mathematics research information

Researchers usually publish their findings in journal articles. Mathematicians are no exception. As a mathematician, you may go directly to the journals to browse and—hopefully—to find relevant information. You may also find interesting references in other journal articles, or you may get reading recommendations from your colleagues. These are all common ways to navigate the realm of your specific topic of research. However, the most efficient and effective way to find information is to search for references in reference databases, also known as bibliographic databases.

A reference database contains information about all journal articles or books published in the major journals and at the most important publishing houses. When you search for a specific topic in a reference database you will find *references* to journal articles or books. A reference to a journal article is simply a statement of the publication details (author, title, name of journal, volume, issue, pages) that allow you to pin-point this specific item and find the full-text. Reference databases generally don't contain the full-text themselves, although often there are links from references in the database to the electronic full text. (A few databases do contain both references and full-text articles; a relevant example is JSTOR, described later in this paper.)

---

[2] Ewing, J. (2002). Twenty Centuries of Mathematics: Digitizing and Disseminating the Past Mathematical Literature. *Notices of the AMS*, 49(7), pp.771-777.

Guthrie, K. (2000). *Revitalizing older published literature: preliminary lessons from the use of JSTOR.* [Online] Available at: http://www.si.umich.edu/PEAK-2000/guthrie.pdf [Accessed Mar 28th, 2009].

The major reference databases in mathematics are *MathSciNet* and *Zentralblatt MATH*, published by the American Mathematical Society and European Mathematical Society respectively. Both these databases are extremely valuable tools when you want to look for references to journal articles in a specific area in mathematics. In the table below you will also find *STMA-Z Statistical Theory*, which is a subset of Zentralblatt MATH in statistics, freely available for all at the time of this writing.

Some African universities already have access to one or both of these major reference databases, but access statistics have shown that the usage at some of these sites is very low or nonexistent. From experience I know that one reason for the low usage is the misconception that the databases contain full-text articles, and, consequently, disappointment when they don't. As explained above, reference databases are useful tools for finding references to relevant articles or books. In some cases, links may be provided to the full text. In cases where such links are not present, just knowing that the publications exist and having the correct references to them is a step forward in your research, and you may be able to get hold of the specific journal articles and books through an interlibrary loan, from colleagues, or through some type of document delivery services such as eJDS, which is described later in this paper.

*By not using reference databases like MathSciNet or Zentralblatt MATH, you may overlook important information and—in the worst case—try to solve problems that have already been worked out by others. If you don't have access to one of these databases, please contact me at anders@golonka.se, and I will find out if you are eligible for free access.*

| Reference database | Address |
| --- | --- |
| MathSciNet | http://www.ams.org/mathscinet/ |
| Zentralblatt MATH | http://www.zentralblatt-math.org/zmath/en/ |
| STMA-Z Statistical Theory | http://www.zentralblatt-math.org/STAT/ |

## "Free" journal resources

When it comes to access to scientific articles, the major divide is between journals paid for through subscriptions and journals that are free. For a researcher at a "northern" university, this distinction normally has no importance. The university library supplies journal articles like water from the tap, and most researchers aren't even aware of which journals are subscribed to and which are free. In the developing world, however, the free journals are extremely important, since subscribing to journals isn't usually feasible because of financial limitations.

The quotation marks on "free" in the header are intentional. "Free" can mean many things. It can mean truly free, that is, free to everyone, everywhere. It can mean available via national provision, through which some publishers provide certain journals free, but only to end-users in low-income countries. Other categories of "free" journals are those that are paid for by donor programmes and by organizations such as the United Nations.

One major force behind the '"truly free" journals is the so-called Open Access movement, which promotes publication of journals that all end-users can access without cost. The term Open Access actually only applies to free-access journal articles provided under a particular license that allows unrestricted derivative use, but the bottom line is that they are free for all. The fact that they are free does not mean they are low quality but rather that they are financed in alternative ways, including through support from institutions or learned societies, author fees, volunteerism or advertising. Another extremely important category of "truly free" journals, in which mathematics is very well covered, is the digitization of old journals, often made possible by governmental support.

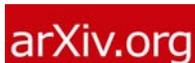

A free resource that is very specific for physics and mathematics is arXiv, which can be described as a database containing working papers or manuscripts for journal articles. The archive was started in 1991 at the Los Alamos National Laboratory in the United States and originally only covered physics. Today you will also find papers there about astronomy, mathematics, computer science, nonlinear science, quantitative biology and statistics.

Authors can upload their scientific papers to arXiv, making them available for free to anyone. In many cases you will find the latest scientific results at arXiv, while you will have to wait for the journal's peer review process before you find the same material in a published form. If you are a mathematician or physicist, you should check arXiv regularly and look for the latest submissions in your topic. If you find a reference to a journal article that you can't get hold of in full text, you may also check arXiv to see if the author perhaps has uploaded the manuscript prior to the real publication.

Time spent exploring the following "true free" resources—Open Access journals, digitized journal archives and electronic preprints—is time well invested.

| Resource | Type | Address |
| --- | --- | --- |
| DOAJ | Open Access journals | http://www.doaj.org/ |
| EMIS | Open Access journals | http://www.emis.de/journals/ |
| DigiZeitschriften | Digitized archive | http://www.digizeitschriften.de/ |
| NUMDAM | Digitized archive | http://www.numdam.org/?lang=en |
| Project Euclid | Digitized archive | http://projecteuclid.org/ |
| Polish Virtual Library of Science | Digitized archive | http://matwbn.icm.edu.pl/index.php?jez=en |
| SwissDML | Digitized archive | http://retro.seals.ch/ |
| arXiv.org | Preprint archive | http://arxiv.org |
| mini-DML | All of the above | http://minidml.mathdoc.fr/ |

## Access methods, authentication, and IP numbers

The "truly free" resources listed in the table above are free to anyone and anywhere. Resources provided by other programmes and initiatives, which are described below, are also free to end-users in all or most African countries. However, there is an importance difference between these two groups of resources. The second group requires some sort of authentication before the user is allowed access.

In other words, in order for you to access a journal that normally requires subscription, the publisher will demand that you prove your identity and eligibility. The two most-common methods for authentication are the use of username/password and IP number control.

The username/password option is a very straightforward method of giving access. All you have to do is supply the correct pair of credentials in order to access the specific resource. This method has its drawbacks, however, since passwords tend to be spread to others who are not eligible for access. When the provider of the specific resource discovers that this is the case, they usually cancel the password and set a new one. The new password must then be distributed to all eligible users, and this may be a slow process.

IP number (or IP address) control is a method in which authentication is based on the computer or computer network from which the researcher connects. All computers on the internet have specific and unique IP numbers. This makes it possible for publishers to determine whether to grant or refuse access to individual users connecting from specific computers by checking the series of numbers supplied by eligible universities and other research institutions.

At smaller universities in Africa, it's common to have only one IP number, the so called public IP number. This public IP number belongs to the gateway computer that connects the university to the outside world. Within the university campus –behind the gateway computer—so-called private network IP numbers are used. These internal numbers are never to be supplied to journals or publishers; only the public IP numbers are of interest.

*To repeat, these ranges of private network IP numbers should never be supplied to publishers for authentication purposes.   They never work!  You have been warned!*

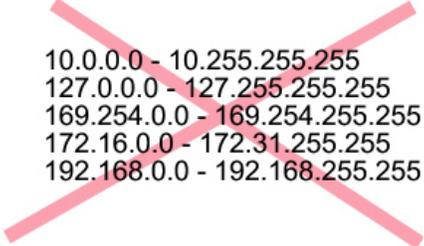

In order for the IP number control system to work smoothly, the public IP number(s) should be fairly stable. In Africa, this is not always the case, since a change of the Internet Service Provider (ISP) also usually means a change of the IP number. African institutions sometimes see an advantage in negotiating terms and prices with a new Internet Service Provider now and then, in order to find a more favorable deal, but this means that the new IP numbers must be supplied to all journals and publishers before access is reestablished.

To complicate this picture a little further, there is a distinction between *static* and *dynamic* IP numbers. In general, there is a world-wide shortage of IP numbers.  In order to cope with this situation, the numbers are sometimes assigned to universities and institutions in a dynamic as opposed to static way. A dynamically assigned IP number may change any time (even though they usually are pretty stable over time). A static number is assigned once and is not supposed to change as long as you have a running contract with an Internet Service Provider, which makes them better for authentication purposes.  The flip-side of the coin is that static numbers are more expensive.

The section above can be summarized as a set of recommendations.  A university should, if possible, stay with one Internet Service Provider and should also try to get static IP numbers.  Please make a note of your institutions' public IP number(s) and make sure it hasn't changed recently.  Never, never supply private network IP numbers (see above) for authentication purposes; they will never work! One economically feasible solution to obtaining internet access with a static public IP number may be for several universities to approach an Internet Service Provider together and negotiate a long-term contract.

## UN sponsored programmes: HINARI, AGORA and OARE

A few of the less-obvious resources for African mathematics are the UN-sponsored access programs HINARI, AGORA and OARE. Although their primary subject coverage is health, agriculture and environmental sciences, they contain surprisingly many high-quality journals in mathematics and statistics.

A search for the subject "Statistics & Mathematics" in HINARI will render you 120 journal titles. A similar search in AGORA on the subjects "Multidisciplinary/ Miscellaneous" gives you some ten journals in the mathematics and statistics area, and a search on "Statistics, Computers & Modeling (Environmental)" in OARE gives you about twenty hits. Especially the Elsevier journals are hard to come by through any other access routes, so do explore these resources with care.

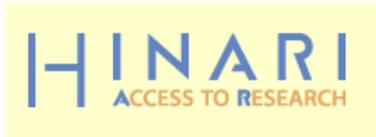

| Resource | HINARI |
|---|---|
| Address | http://www.who.int/hinari/ |
| Eligibility and costs | Country specific.<br>GNI per capita less than USD 1250 = free access<br>GNI per capita between USD 1250 and 3500 = USD 1000/year |
| Type of material | Electronic journals; health sciences *but also many mathematics journals* |
| Access control | Username and password |

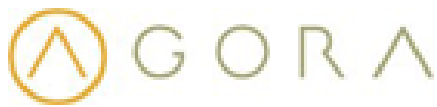

| Resource | AGORA |
|---|---|
| Address | http://www.aginternetwork.org/en/ |
| Eligibility and costs | Country specific.<br>GNI per capita less than USD 1000 = free access<br>GNI per capita between USD 1000 and 3000 = USD 1000/year |
| Type of material | Electronic journals; agricultural sciences *but also some statistics and probability journals* |
| Access control | Username and password |

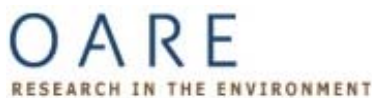

| Resource | OARE |
|---|---|
| Address | http://www.oaresciences.org/en/ |
| Eligibility and costs | Country specific.<br>GNI per capita less than USD 1250 = free access<br>GNI per capita between USD 1250 and 3500 = USD 1000/year |
| Type of material | Electronic journals; environmental sciences *but also some statistics and probability journals* |
| Access control | Username and password |

Most developing country institutions already have usernames and passwords for these UN access programmes, but since mathematicians are not the major target group, the information may not have reached you. In order to get the valid usernames and passwords for HINARI, AGORA and OARE, please contact your local university library. In the unlikely event that your university isn't registered, encourage your main library to do so.  If push comes to shove (nothing else works), register your mathematics department separately. Information on how to register is available on the web pages of the respective programme, under a link called "Register" in the left hand menu.

*Please don't give up! The mathematics journals in these UN programmes are too good to be overlooked or taken lightly!*

## Journals available via donor programmes or national provision

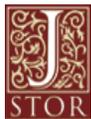

JSTOR is an archive of scanned journals in a number of different disciplines, including mathematics. The cost for this resource is normally quite high, but in 2006 JSTOR instated a special scheme waiving the charges for all academic and not-for-profit institutions on the African continent.

| **Resource** | **JSTOR** |
|---|---|
| **Address** | http://www.jstor.org/ |
| **Eligibility and costs** | Not-for-profit institutions on the African continent = free access |
| **Type of material** | Electronic journals; most disciplines including a *notably good mathematics and statistics collection* |
| **Access control** | IP control |

JSTOR may be the single most valuable resource for African mathematicians, since the titles and coverage of the journals in the "Mathematics & Statistics" collection is impressive. Worth mentioning are all major journals from the American Mathematical Society, the American Statistical Society and the Institute of Mathematical Statistics from the publication start until about 2005.  Note that JSTOR is an archive, with articles added a few years after their publication.

Access control to JSTOR is based on IP number, and all institutions have to apply,[3] first by filling in the "JSTOR Network Verification Form", and then submitting a printed contract in two copies. Apart from the mathematics and statistics journals, JSTOR contains publications in many other areas.  So please don't miss the opportunity to sign up for this resource. And, once access is established, do inform your colleagues in other departments at your university.

---

[3] http://www.jstor.org/page/info/participate/new/fees/africanAccess.jsp

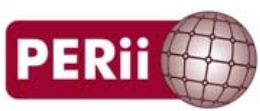

INASP is a non-governmental organization based in Oxford, UK, which handles access to scientific publications for developing countries in a very broad sense. INASP deals with training, advising, ICT-issues and publishing. Also—most important in the context of this paper—it negotiates fair prices for journals and databases from many of the worlds major publishing houses via the INASP/PERii project.

| Resource | INASP PERii |
|---|---|
| Address | http://www.inasp.info/ |
| Eligibility | Country specific |
| Type of material | Journals, all subjects including mathematics and statistics |
| Access control | IP control |

Using PERii resources involves both registration with PERii and with the provider or publisher of each specific resource. In order to find out what is available for a specific country, go to the INASP webpages at http://inasp.info and choose a country in the "Country Finder" drop-down menu. You will then see a list of resources that are available for your country. Click on a resource to obtain more information on how to proceed to get access (usually under the header "Registration Information"). By clicking the link to "INASP registration system" you will reach a list of eligible countries. Click again on your country, and a list of institutions registered with PERii will appear. If your institution is in the list you can continue by clicking on the link and follow the instructions for registering for the resource. If your institution doesn't appear in the list you may add your details by clicking "add your institution" further down the page.

Registration for INASP PERii resources is something which ideally should be done by your university library. As in the case of the UN sponsored programmes above; if you fail to get support from your local library, you may register as a mathematics department.

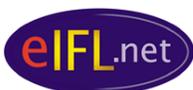

eIFL is an not-for-profit organization which assists lower- and middle-income countries in forming national library consortia and negotiating affordable subscription prices for electronic journals. The following African countries have formed library consortia with the help from eIFL.net: Botswana, Cameroon, Egypt, Ethiopia, Ghana, Kenya, Lesotho, Malawi, Mali, Mozambique, Nigeria, Senegal, South Africa, Sudan, Swaziland, Zambia and Zimbabwe. If you live in one of the countries above, you can find the eIFL.net local coordinator by going to http://www.eifl.net/cps/sections/country .

| Resource | eIFL.net |
|---|---|
| Address | http://eifl.net |
| Eligibility | Country specific, see http://www.eifl.net/cps/sections/country |
| Type of service | Library consortium for electronic journals |
| Access control | IP control |

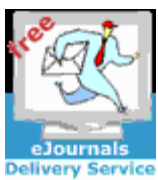 **The electronic Journals Delivery Service (eJDS)**

eJDS is a programme initiated by the Abdus Salam International Centre for Theoretical Physics, Trieste, Italy. Developing-country physicists and mathematicians can register for access to a number of journals from publishers such as the American Mathematical Society, Institute of Physics, American Physical Society, American Institute of Physics, SPIE and Optical Society of America. After a successful registration, researchers can log on to the service with their email address, and the papers requested are sent back to the researchers to the same email address. The advantage of sending papers via email is that the bandwidth required is minimal and that the service can be used anywhere, including at the local internet café. There is, however, a limit on how many papers the individual researcher can download; 3 per day, 12 per month and 100 per year.

| Resource | eJDS |
|---|---|
| Address | http://ejds.org/ |
| Eligibility | Country specific |
| Type of material | Journals in physics and mathematics |
| Access control | Registration required, delivery of articles via email |

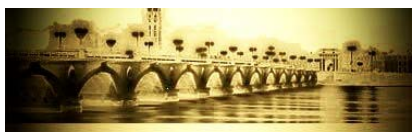 **Bordeaux*thèque***

Another valuable free document delivery resource - very similar to the eJDS above - is a service called **Bordeaux*thèque***, run by l'Université Bordeaux 1 in France. Registration is required, and once accepted, the developing country researcher will be assigned a personal member ID. Ordering of papers can then be done through an online form. As in the case of eJDS, documents are delivered via email. Bordeauxthèque can at the time of this writing supply documents from 287 current and 180 discontinued major journals in mathematics and computer science.

It is worth mentioning that there is a considerable overlap among the journals provided by different providers like HINARI, JSTOR, eJDS, Bordeauxthèque and others. So, if you don't find the specific volume in one of the resources, you should also check the others.

| Resource | Bordeaux*thèque* |
|---|---|
| Address | http://bordeauxtheque.math.u-bordeaux1.fr/ |
| Eligibility | Researchers from public and non-commercial institutions in developing countries. |
| Type of material | Journals in mathematics and computer science |
| Access control | Registration required, delivery of articles via email |

**Learning more**

In this paper, I have briefly introduced a selection of the many valuable information resources available free-of-charge for mathematicians in Africa. If you want to learn more, I recommend that you do two things. First, go online to my website, "e-Math for Africa", at http://math.golonka.se and go through the links provided there.

Second, contact me about INFORM, the International Network for Online Resources and Materials, a programme based at Uppsala University in Sweden. Although I primarily work as a librarian, I also am employed part-time as an INFORM trainer. Martha Garrett (the INFORM director) and I have carried out information training for mathematicians in Africa and have co-authored a compendium on accessing mathematics and physics information resources in Africa. If you send me an email at anders@golonka.se, I can provide you with more information about training possibilities and about accessing the latest version of the compendium, which is both cost-free and copyright-free.

And remember, online mathematics information in Africa is not served on a silver platter. You will have to work hard for your treat. Don't give up!